\documentclass[10pt]{article}
\usepackage{amsmath,amssymb,amsthm,mathtools}
\usepackage[a4paper,margin=1in]{geometry}
\usepackage{hyperref}
\usepackage{authblk}
\newtheorem{theorem}{Theorem}[section]
\newtheorem{lemma}[theorem]{Lemma}
\newtheorem{proposition}[theorem]{Proposition}
\newtheorem{corollary}[theorem]{Corollary}
\theoremstyle{remark}
\newtheorem{remark}[theorem]{Remark}

\newcommand{\E}{\mathbb E}
\newcommand{\Pp}{\mathbb P}
\newcommand{\Lop}{\mathcal L}
\newcommand{\Qop}{\mathcal Q}
\newcommand{\one}{\mathbf 1}
\newcommand{\norm}[1]{\lVert #1\rVert}
\newcommand{\Dcone}{\mathcal D}

\usepackage{cite}

\begin{document}
	
	\title{Exponential loss of memory for random Lasota--Yorke convex maps with local average expansion}
	
	\author{
		Hongfei Cui\thanks{School of Mathematics and Systems Science, Hubei Province Key Laboratory of Systems Science in Metallurgical Process, Wuhan University of Science and Technology, Wuhan 430065, People's Republic of China \\
			Email: cui@wipm.ac.cn}
		\qquad
		Juho  Lepp\"{a}nen\thanks{Department of Mathematics, Tokai University, Kanagawa 259-1292, Japan \\
			Email: leppanen.juho.heikki.g@tokai.ac.jp}
	}
	
	\date{\today}
	
	\maketitle

	\begin{abstract}
		We study i.i.d. random compositions of piecewise convex interval maps whose first branch is full.
		Uniform expansion is not required, and the maps may have contracting branches.
		Under a local average expansion condition at the origin, we prove both quenched and annealed exponential loss of memory for initial densities of bounded variation. The proof is based on a coupling argument. As applications, we obtain the quenched almost sure invariance principle and a rate of convergence in the quenched central limit theorem.
	\end{abstract}
	

	\section{Introduction}
	
	Statistical loss of memory in a dynamical system refers to
	different initial probability distributions becoming
	asymptotically indistinguishable as they evolve under the dynamics.
	More precisely, for a sequence of maps with associated transfer
	operators $(L_n)_{n\geq 1}$, we consider loss of memory in the sense
	that
	\begin{equation}\label{eq:ml}
		\|L_n\cdots L_1(f-g)\|_{L^1(m)} \longrightarrow 0
		\qquad\text{as } n\to\infty
	\end{equation}
	for probability densities $f$ and $g$ in a class of suitably regular functions,
	where $m$ is a reference measure on the state space,
	such as Lebesgue measure. This notion applies both to iterates of a single map
	and to nonautonomous compositions, in which the map changes with
	time. In the latter setting, the evolving densities need not
	converge to a fixed equilibrium, but may nevertheless approach
	one another. The rate of convergence in \eqref{eq:ml} measures how quickly
	the initial distribution is forgotten. Estimates of this rate
	are an important ingredient in proofs of various statistical limit
	theorems, including central limit theorems and invariance
	principles.

	For nonautonomous compositions of uniformly expanding maps,
	exponential rates of memory loss can often be established using
	spectral methods or contraction in projective metrics.
	These approaches have been developed in \cite{GOT13, CR07,B99}, among others,
	and yield exponential convergence in norms stronger
	than $L^1$, such as the bounded variation norm or quasi-H\"older
	norms.
	
	Another approach is coupling, which directly estimates the
	$L^1$-distance between evolving densities by successively
	removing common mass. Once a common part of two densities has
	been identified, it has the same subsequent evolution and
	therefore makes no contribution to their difference.
	The uncoupled remainders are then evolved until the procedure
	can be repeated. Coupling has been used to establish rates
	of memory loss for nonautonomous systems in both uniformly
	expanding or hyperbolic \cite{SYZ13, OSY09} and nonuniformly expanding
	\cite{KL21, BBR19}
	settings.

	For sequential compositions of Lasota--Yorke convex maps on the unit interval, Cui
	\cite{Cui2021} established exponential loss of memory in $L^1$
	by combining coupling with preservation of non-increasing
	densities and suitable cones by the transfer operators.
	Large deviations and central limit theorems for sequential
	and random compositions
	were obtained in \cite{CS24, CW25}.
	These maps are
	piecewise increasing and convex, and their first branch is full,
	mapping its domain onto the whole interval. The
	derivative at the origin of the first branch is assumed to satisfy a uniform lower
	bound strictly greater than one. By convexity, this gives uniform
	expansion throughout that branch, while the other branches need
	not be expanding.
	
	In this paper, we study i.i.d. compositions of Lasota--Yorke convex
	maps and replace the uniform expansion condition on the first branch in
	\cite{Cui2021} by an expansion-on-average
	condition.
	More precisely, let
	\[
	\lambda_\xi=T'_{\xi,1}(0+)
	\]
	be the derivative at the origin of the first, full branch of the random
	map $T_\xi$. Our main assumption is
	\begin{align}\label{eq:intro_local}
		\theta:=\E(\lambda_\xi^{-1})<1.
	\end{align}
	We refer to this as \emph{local average expansion}. The condition is local
	because it concerns only the derivative of the first branch at the origin.
	Individual maps may therefore be contracting near the origin,
	and the remaining branches may also be contracting. In particular,
	we do not impose an average expansion condition on the minimum rate of expansion
	over all branches.
	
	Under our standing assumptions,
	we prove both quenched and annealed exponential loss of memory in $L^1$. The annealed
	estimate controls the expected $L^1$-distance between two evolving
	densities, whereas the quenched estimate holds along almost every
	realization of the random environment, with a deterministic exponential rate and an integrable
	random coefficient. We first prove these estimates for
	bounded
	non-increasing probability densities and then
	extend them to densities with bounded variation.
	As a consequence of the annealed estimate we establish the
	existence and uniqueness
	of the invariant probability density for
	the averaged transfer operator. As applications of the quenched exponential memory loss,
	we show the existence of the sample measures, the
	quenched almost sure invariance principle and a rate of convergence in the quenched
	central limit theorem. For the latter two results we apply abstract theorems
	in \cite{DHL2026}.
	
	We emphasize that quenched limit theorems for expanding-on-average
	and nonuniformly expanding random dynamical systems have been studied
	extensively; see, for example,
	\cite{DS23, DHS23, H23, DH21}.
	In examples considered in \cite{DS23, DHS23},
	the average expansion condition is formulated
	in terms of the minimum expansion rate over the phase space.
	For the family of convex maps considered here, we obtain quenched limit theorems
	under the local condition
	\eqref{eq:intro_local}, without imposing average expansion
	on the remaining branches.

	Our proof adapts the coupling approach of Stenlund and Sulku
	\cite{StenlundSulku2014}, who considered random compositions of
	smooth circle maps expanding on average. Their expansion hypothesis
	is global in space: it is formulated in terms of
	$\inf_x |T_\xi'(x)|$, together with a moment condition on distortion.
	We instead use the convexity of the maps and the preservation of
	non-increasing densities by the transfer operator, together with the local
	expansion on average condition \eqref{eq:intro_local}.
	The regularity quantity we analyze is the \emph{height} of a
	non-increasing density $f$, namely
	\[
	f(0+)=\operatorname{ess\,sup}f.
	\]
	Controlling this quantity under iteration of the transfer operators
	prevents the mass from becoming too concentrated near the origin
	and makes it possible to couple a common part of two densities
	after evolving them along a suitable block of maps.
	
	There are two difficulties in implementing this procedure.
	First, an unfavorable sequence of maps can increase the height
	of a density, and removing common mass can further increase the
	height of the normalized uncoupled remainder. It is therefore
	necessary to control the remaining mass and its concentration near the origin
	together to sustain the coupling procedure. Second, the lower
	bound used for coupling depends on the block of maps over which
	the densities evolve. If coupling is attempted after a variable
	waiting time, the available common fraction need not have a
	uniform positive lower bound.
	
	To address these difficulties, we perform coupling attempts at
	intervals of a fixed block length.
	A block is called good when the contraction of the inverse
	of the first branch over that block satisfies a prescribed bound.
	The condition $\theta<1$ and Markov's inequality,
	together with independence, imply that each block is good with
	a conditional probability bounded below by a positive constant.
	If the normalized uncoupled densities have controlled height,
	a good block produces a deterministic lower bound for the
	evolved densities, allowing a fixed proportion of their remaining
	mass to be coupled. If the height is too large, no coupling is
	performed, but the height decreases in conditional expectation.
	
	A weighted quantity involving the uncoupled mass and the height
	of the normalized remainder therefore contracts uniformly from
	one block to the next in conditional expectation. This yields
	both annealed exponential loss of memory and a quenched
	exponential estimate with an integrable random coefficient
	for non-increasing densities. By expressing densities of bounded
	variation as differences of non-increasing functions, we extend
	these estimates to initial densities with bounded variation.
	
	The paper is organized as follows. Section~2 introduces the random
	Lasota--Yorke convex family, formulates the local average expansion
	condition, and states the main theorem. Section~3 develops estimates for densities evolving under the
	random dynamics, including a lower bound over blocks of fixed
	length, and formulates the coupling procedure.
	Section~4 proves a conditional contraction estimate over a
	single block and derives annealed and quenched exponential
	loss of memory. Section~5 treats bounded-variation initial densities
	and the averaged transfer operator. Section~6 establishes the
	existence of sample measures and gives applications to quenched
	limit theorems.

	\section{Setting and main result}
	
	Let $(\Omega,\mathcal A,\eta)$ be a probability space and let
	$
	\Pp=\eta^{\otimes\mathbb N}
	$
	be the product law of the i.i.d.\ environment coordinates. Write
	$
	\omega=(\omega_1,\omega_2,\ldots).
	$
	Let $T_\xi:[0,1]\to[0,1]$, $\xi\in\Omega$, be a random family. We
	assume that $(\xi,x)\mapsto T_\xi(x)$ is jointly measurable. The following assumptions hold
	for $\eta$-almost every $\xi$.
	
	\medskip
	\noindent\emph{(H1) Piecewise convexity.}
	There is a finite partition
	$
	0=a_{\xi,0}<a_{\xi,1}<\cdots<a_{\xi,N_\xi}=1
	$
	such that the restriction of $T_\xi$ to each partition interval is
	$C^1$ and convex. Denote $T_{\xi, i}:=T_{\xi}|_{(a_{\xi,i-1}, a_{\xi,i})}$ for $i=1, \dots, N_{\xi}$.
	The first branch is full:
	\[
	T_{\xi,1}(0+)=\lim_{x\downarrow0}T_{\xi,1}(x)=0,
	\qquad
	T_{\xi,1}\bigl((0,a_{\xi,1})\bigr)=(0,1).
	\]
	For every $i\ge2$,
	\[
	\lim_{x\downarrow a_{\xi,i-1}}T_{\xi,i}(x)=0,
	\qquad
	T'_{\xi,i}(a_{\xi,i-1}+)>0.
	\]
	Consequently, convexity and positivity of the left derivative make
	every branch strictly increasing on its domain.
	
	\medskip
	\noindent\emph{(H2) Uniform bounds.}
	There are constants $a\in(0,1)$, $M_1>1$, and $M>0$, independent of
	$\xi$, such that
	\[
	a_{\xi,1}\ge a,
	\qquad
	T'_{\xi,1}(a_{\xi,1}-)\le M_1,
	\]
	and
	\[
	\sum_{i=2}^{N_\xi}
	\frac{1}{
		T'_{\xi,i}(a_{\xi,i-1}+)\,a_{\xi,i-1}}
	\le M.
	\]

	\medskip
	\noindent\emph{(H3) Local average expansion.}
	Put
	$
	\lambda_\xi=T'_{\xi,1}(0+)\in(0,\infty),\ 
	\theta=\int_\Omega\lambda_\xi^{-1}\,d\eta(\xi).
	$
	We assume that
	\begin{equation}\label{eq:theta}
		\theta<1.
	\end{equation}

	Let $L_\xi$ denote the transfer operator of $T_\xi$ with
	respect to Lebesgue measure $m$. Thus
	\[
	\int_0^1(L_\xi f)g\,dm
	=
	\int_0^1f(g\circ T_\xi)\,dm
	\]
	for $f\in L^1(m)$ and $g\in L^\infty(m)$.
	For
	$
	\omega=(\omega_1,\omega_2,\ldots) \in \Omega^{\mathbb N},
	$
	write
	\[
	\Lop_\omega^n
	=
	L_{\omega_n}\cdots L_{\omega_1},
	\qquad n\ge1,
	\]
	and let $\Lop_\omega^0$ be the identity. Similarly, write
	\[
	T_\omega^n
	=
	T_{\omega_n}\circ\cdots\circ T_{\omega_1}.
	\]
	
	All function equalities are understood modulo Lebesgue null sets.
	As usual, elements of $L^1(m)$ are identified with their equivalence classes
	modulo equality almost everywhere.
	We say that $f\in L^1(m)$ is non-increasing if it admits a
	non-increasing representative. For each such class, we choose its
	unique right-continuous non-increasing representative on $(0,1)$.
	Recall that the set of non-increasing elements of $L^1(m)$ is
	closed in $L^1(m)$.
	
	For non-negative, non-increasing $f\in L^1(m)$, we write
	$$
	f(0+)=\mathop{\rm ess\,sup} f.
	$$
	We use $\E$ to denote the expectation with respect to $\Pp$.
	We will use the bounded variation norm
	\[
	\|\psi\|_{\mathrm{BV}}
	=
	\|\psi\|_{L^\infty}+\operatorname{Var}(\psi)
	\]
	and  the Lipschitz norm
	\[
	\|g\|_{\mathrm{Lip}}
	=
	\|g\|_\infty+\operatorname{Lip}(g).
	\]
	Here,
	\begin{align*}
		\operatorname{Var}(\psi) = \inf_{ \text{ $\bar{\psi} = \psi$ $m$-a.e.}  }
		\sup \biggl\{
		\sum_{j=0}^{t-1} | \bar{\psi}( x_{j+1} ) -  \bar{\psi}( x_{j} )  | \: : \:
		0 = x_0 < \ldots < x_t = 1
		\biggr\}
	\end{align*}
	and
	\begin{align*}
		\operatorname{Lip}(g) = \sup_{ \substack{ x \neq y \\ x,y \in [0,1] } }
		\frac{ |g(x) - g(y)| }{|x-y|}.
	\end{align*}
	For a non-negative and non-increasing $f \in L^1(m)$ we have
	$
	\operatorname{Var}(f) \le f(0+).
	$
	We define
	$$
	BV([0,1]) = \{  f \in L^\infty \: : \: \text{Var}(f) < \infty \}.
	$$
	
	For $B\ge1$, define
	\[
	\Dcone_B
	=
	\left\{
	f\ge0:
	f\text{ is non-increasing},\
	\int_0^1f\,dm=1,\
	f(0+)\le B
	\right\}.
	\]
	The parameter $B$ controls the height of a probability density at the
	origin. Although the transfer operators preserve non-negativity,
	integrals, and monotonicity, a fixed cone $\Dcone_B$ need not be
	preserved by every $L_{\xi}$, in particular because $T_{\xi}$ may
	be locally contracting at the origin.
	Controlling this possible growth of
	the height under the action of the transfer operators
	will be an essential part of our  coupling argument.
	
	\begin{theorem}
		\label{thm:main}
		Under assumptions \emph{(H1)--(H3)} and the preceding measurability and
		i.i.d.\ assumptions, there exists $\rho\in(0,1)$, depending only on
		$a,M_1,M,\theta$, such that for every $B\ge1$ there is $C_B<\infty$
		satisfying
		\[
		\E\left[
		\norm{\Lop_\omega^n(f-g)}_{L^1}
		\right]
		\le C_B\rho^n,
		\qquad n\ge0,
		\]
		for every $f,g\in\Dcone_B$.
		
		Moreover, there exists $c>0$ such that for every $B\ge1$
		there is a non-negative random variable
		\[
		C_{B,c}\in L^1(\Pp)
		\]
		such that, almost surely,
		\[
		\norm{\Lop_\omega^n(f-g)}_{L^1}
		\le C_{B,c}(\omega)e^{-cn},
		\qquad n\ge0,
		\]
		for all $f,g\in\Dcone_B$. The constant $c$ is independent of $B$.
	\end{theorem}
	
	The first estimate in Theorem~\ref{thm:main} is the annealed
	exponential loss of memory. The second is the quenched estimate: for
	almost every fixed $\omega$, the memory of the initial density is
	lost exponentially, with a random but integrable coefficient.

	\section{The coupling argument}
	We now adapt the coupling construction of~\cite[Section~5.2]{StenlundSulku2014} to the present random piecewise convex
	interval maps with local average expansion.  Throughout the proof, $\omega$ is reserved for the
	environment, while $f,g,$ denote densities.
	The argument has three ingredients: preservation of non-increasing
	densities by the transfer operators, a fixed-block lower bound,
	and a recursion for the mass and height of the uncoupled remainders.
	
	\subsection{Cone and height estimates}
	For a single map $T_\xi$, put $I_{\xi,i}=(a_{\xi,i-1},a_{\xi,i})$ and $J_{\xi,i}=T_\xi(I_{\xi,i})$. Positivity of the left derivative and convexity make every branch strictly increasing. Every $J_{\xi,i}$ is an interval starting at zero, and $J_{\xi,1}=(0,1)$. With $g_{\xi,i}:J_{\xi,i}\to I_{\xi,i}$ its inverse,
	\[
	L_\xi f(y)=\sum_{i:y\in J_{\xi,i}}
	\frac{f(g_{\xi,i}(y))}{T'_{\xi,i}(g_{\xi,i}(y))}\quad\text{a.e.}
	\]
	Each summand is non-increasing for non-negative non-increasing $f$: both factors $f\circ g_{\xi,i}$ and $1/(T'_{\xi,i}\circ g_{\xi,i})$ are non-increasing. Its zero extension beyond $J_{\xi,i}$ is also non-increasing. Thus $L_\xi$ preserves the non-negative decreasing cone. Change of variables on each branch gives $\int L_\xi f\,dm=\int f\,dm$ and positivity yields $\|L_\xi h\|_1\le\|h\|_1$.
	\begin{lemma}\label{lem:onestep}
		If $f\ge0$ is non-increasing, then
		\[
		(L_\xi f)(0+)
		\le \lambda_\xi^{-1}f(0+)+M\int_0^1f\,dm.
		\]
	\end{lemma}
	\begin{proof}
		For $x>0$, monotonicity gives
		\[
		f(x)\le x^{-1}\int_0^1f\,dm.
		\]
		The contribution of the first branch to $\lambda_\xi^{-1}f(0+)$ is bounded by
		$\lambda_\xi^{-1}f(0+)$. For $i\ge2$, the contribution of the
		$i$th branch to $(L_\xi f)(0+)$ is bounded by
		\[
		\frac{f(a_{\xi,i-1}+)}{T'_{\xi,i}(a_{\xi,i-1}+)}
		\le
		\frac{\int f\,dm}{T'_{\xi,i}(a_{\xi,i-1}+)a_{\xi,i-1}}.
		\]
		Summation combined with (H2) gives the result.
	\end{proof}
	Given $\omega=(\omega_1, \omega_2, \dots) \in \Omega^{\mathbb{N}}$, for $n\ge1$, put
	\[
	S_n(\omega)=\prod_{j=1}^n\lambda_{\omega_j}^{-1},
	\qquad
	R_n(\omega)=\sum_{i=1}^n\prod_{j=i+1}^n\lambda_{\omega_j}^{-1}.
	\]
	Iteration of Lemma~\ref{lem:onestep} yields
	\begin{equation}\label{eq:height}
		(\Lop_\omega^nf)(0+)
		\le S_n(\omega)f(0+)+MR_n(\omega)\int f\,dm.
	\end{equation}
	
	The height estimate \eqref{eq:height} controls how concentrated a decreasing density may
	become at the origin
	under a composition, but the coupling step requires a uniform
	lower bound.  The next lemma provides precisely this complementary
	estimate.  The resulting constant
	$\kappa_\ell$ is independent of the density and of the particular
	$\omega$, and will be the fixed fraction removed at each
	successful coupling.
	
	\begin{lemma}\label{lem:lower}
		Let $A\ge1$ and set $t_A=a/(2A)$. If $f\in\Dcone_A$ and a block $\omega$ of length $\ell$ satisfies
		\[
		S_\ell(\omega)\le t_A,
		\]
		then
		\[
		L_{\omega_\ell}\cdots L_{\omega_1}f\ge \kappa_\ell
		\quad\text{a.e.},
		\]
		where
		\[
		\kappa_\ell=\min\left\{\frac12,
		M_1^{-\ell}\frac{1-a}{1-t_A}\right\}.
		\]
	\end{lemma}
	\begin{proof}
		Since $f(0+)\le A$ and $f$ is non-increasing,
		\[
		1=\int_0^{t_A}f+\int_{t_A}^1f
		\le At_A+(1-t_A)f(t_A).
		\]
		Thus
		\[
		f(t_A)\ge\frac{1-a}{1-t_A}.
		\]
		For $y\in(0,1)$ let $b_j=g_{\omega_j,1}$ and set
		\[
		z_\ell=y,\qquad z_{j-1}=b_j(z_j)\ (j=\ell,\ldots,1),
		\qquad x_y=z_0=b_1\circ\cdots\circ b_\ell(y).
		\]
		All these inverse branches exist because the first branch of each map is full.
		Convexity, the identity $T_{\omega_j,1}(0+)=0$, and the
		derivative bound $T'_{\omega_j,1}\ge\lambda_{\omega_j}$
		imply $b_j(z)\le\lambda_{\omega_j}^{-1}z$, so
		\[
		x_y\le S_\ell(\omega) y\le t_A.
		\]
		The intermediate point $z_{j-1}$ lies in the first branch of $T_{\omega_j}$. Hence
		\[
		(T_\omega^\ell)'(x_y)
		=\prod_{j=1}^\ell T'_{\omega_j,1}(z_{j-1})\le M_1^\ell.
		\]
		Keeping only this branch gives
		$\Lop_\omega^\ell f(y)\ge M_1^{-\ell}f(x_y)
		\ge M_1^{-\ell}(1-a)/(1-t_A)$. This proves the asserted a.e. bound.
	\end{proof}
	\subsection{The fixed-block coupling step}

	Set
	\[
	D=\frac{M}{1-\theta},
	\qquad
	A=\max\{2,8D\}.
	\]
	Choose $\ell$ so large that
	\[
	q:=\theta^\ell\le\min\left\{\frac14,\frac{a}{4A}\right\}.
	\]
	Fix this value of $\ell$ for the remainder of the proof.
	For $\omega=(\omega_1, \omega_2,\dots)$ and $k\geq 1$, we consider the first $k\ell$ terms $$(\omega_1, \dots, \omega_{\ell},  \omega_{\ell+1}, \dots, \omega_{2\ell}, \dots, \omega_{(k-1)\ell+1}, \dots, \omega_{k\ell}),$$
	which we divide into $k$ blocks with every
	block consisting of exactly $\ell$ consecutive maps.
	Let
	\[
	\mathcal F_k=\sigma(\omega_1,\ldots,\omega_{k\ell})
	\]
	be the block filtration. Here, by a slight abuse of notation,
	$\omega_j$ also denotes the coordinate projection
	$\omega\mapsto\omega_j$ on $\Omega^{\mathbb N}$.
	For block $k$, let
	\[
	\widehat S_k
	=\prod_{j=(k-1)\ell+1}^{k\ell}\lambda_{\omega_j}^{-1},
	\]
	\[
	\widehat R_k=\sum_{i=(k-1)\ell+1}^{k\ell}
	\prod_{j=i+1}^{k\ell}\lambda_{\omega_j}^{-1}.
	\]
	We call the $k$-th block $(\omega_{(k-1)\ell+1}, \dots, \omega_{k\ell})$ is {\it good} if $\widehat S_k\le \frac{a}{2A}$.
	We define the good-block
	\[
	G_k=\left\{\widehat S_k\le \frac{a}{2A}\right\}.
	\]

	In the construction of \cite[Section 5.2]{StenlundSulku2014}, a coupling condition first puts the
	two evolved densities in a regularity class with a uniform positive
	lower bound, one then removes the same constant part from both densities.
	Here the role of the regularity parameter is played by height divided by
	uncoupled mass.  Since that ratio need not remain bounded after a
	coupling, we record both quantities instead of asserting that every
	normalized remainder returns immediately to a fixed cone.

	We next describe exactly what is coupled.  Start with two probability
	densities $f,g\in\Dcone_B$.  At the end of block $k$, corresponding to
	the ordinary iteration time $k\ell$, the functions $u_k$ and $v_k$
	represent the parts of the two evolved densities that have not yet
	been matched.  We denote by $m_k$ {\it the mass} of $u_k$ or $v_k$. Denote by $h_k$
	a common upper bound for  $u_k$ and  $v_k$  at the origin:
	\begin{equation}\label{eq:coupling-invariants}
		u_k,v_k\ge0\ \text{and are non-increasing},\qquad
		\int u_k\,dm=\int v_k\,dm=m_k,
		\qquad
		u_k(0+),v_k(0+)\le h_k.
	\end{equation}
	The initial values are
	\[
	u_0=f,\qquad v_0=g,\qquad m_0=1,\qquad h_0=B.
	\]
	Thus no mass has been coupled at time zero.
	We now construct these objects recursively.
	\paragraph{The first block: $k=1$.}
	The first block operator is
	\[
	\mathcal P_1=L_{\omega_\ell}\cdots L_{\omega_1}.
	\]
	Since $m_0=1$ and $h_0=B$, the coupling condition at the first block is
	\begin{equation}\label{eq:first-coupling-condition}
		\mathcal C_1(A)
		=\{h_0\le Am_0\}\cap G_1
		=\{B\le A\}\cap G_1,
		\qquad I_1=\one_{\mathcal C_1(A)}.
	\end{equation}
	Notice that a good first block is not by itself sufficient: if
	$B>A$, then $I_1=0$ and no coupling is attempted at this block.
	
	Define
	\begin{equation}\label{eq:first-coupling-recursion}
		\begin{aligned}
			u_1&=\mathcal P_1f-\kappa_\ell I_1,\\
			v_1&=\mathcal P_1g-\kappa_\ell I_1,\\
			m_1&=1-\kappa_\ell I_1,\\
			h_1&=\widehat S_1B+M\widehat R_1.
		\end{aligned}
	\end{equation}
	
	We check both possible values of $I_1$.  If $I_1=0$, then
	\[
	u_1=\mathcal P_1f,\qquad v_1=\mathcal P_1g,
	\qquad m_1=1.
	\]
	The transfer operator is positive, preserves integrals, and preserves
	non-increasing functions.  Hence $u_1,v_1$ are non-negative and
	non-increasing, and both have integral $m_1=1$.  Moreover, the block
	height estimate \eqref{eq:height} gives
	\[
	u_1(0+)=(\mathcal P_1f)(0+)
	\le \widehat S_1 f(0+)+M\widehat R_1\int f\,dm
	\le \widehat S_1B+M\widehat R_1=h_1,
	\]
	and the same argument gives $v_1(0+)\le h_1$.
	
	Suppose instead that $I_1=1$.  Then $B\le A$, so
	$f,g\in\Dcone_A$, and $G_1$ occurs.  Lemma~\ref{lem:lower} therefore
	gives
	\begin{equation}\label{eq:first-common-lower-bound}
		\mathcal P_1f\ge\kappa_\ell,
		\qquad
		\mathcal P_1g\ge\kappa_\ell.
	\end{equation}
	It follows that $u_1,v_1\ge0$.  They are non-increasing because a
	constant has been subtracted from two non-increasing functions.  Since
	$m([0,1])=1$ and the transfer operator preserves integrals,
	\[
	\int u_1\,dm
	=\int\mathcal P_1f\,dm-\kappa_\ell
	=1-\kappa_\ell=m_1,
	\]
	and likewise $\int v_1\,dm=m_1$.  Finally,
	\[
	u_1(0+)=(\mathcal P_1f)(0+)-\kappa_\ell
	\le(\mathcal P_1f)(0+)\le h_1,
	\]
	and similarly $v_1(0+)\le h_1$.  Thus all four assertions in
	\eqref{eq:coupling-invariants} hold at $k=1$.  Notice also that
	$m_1\ge1-\kappa_\ell>0$, so we can normalize at the next block.
	
	\paragraph{Induction step: from $k-1$ to $k$.}
	Fix $k\ge2$.  Assume that the objects
	\[
	u_s,v_s,m_s,h_s,\qquad 0\le s\le k-1,
	\]
	have already been defined, and assume in particular that
	\eqref{eq:coupling-invariants} holds at time $k-1$ and that
	$m_{k-1}>0$.  Define the next block operator by
	\[
	\mathcal P_k
	=L_{\omega_{k\ell}}\cdots L_{\omega_{(k-1)\ell+1}}.
	\]
	The coupling condition for this block is
	\begin{equation}\label{eq:coupling-condition}
		\mathcal C_k(A)
		=\{h_{k-1}\le Am_{k-1}\}\cap G_k,
		\qquad I_k=\one_{\mathcal C_k(A)}.
	\end{equation}
	We then define
	\begin{equation}\label{eq:coupling-recursion}
		\begin{aligned}
			u_k&=\mathcal P_ku_{k-1}-\kappa_\ell m_{k-1}I_k,\\
			v_k&=\mathcal P_kv_{k-1}-\kappa_\ell m_{k-1}I_k,\\
			m_k&=m_{k-1}(1-\kappa_\ell I_k),\\
			h_k&=\widehat S_kh_{k-1}+M\widehat R_km_{k-1}.
		\end{aligned}
	\end{equation}
	
	We prove that these definitions preserve the induction hypotheses.
	First consider $I_k=0$.  No common part is removed, so
	\[
	u_k=\mathcal P_ku_{k-1},\qquad
	v_k=\mathcal P_kv_{k-1},\qquad
	m_k=m_{k-1}>0.
	\]
	Positivity and preservation of monotonicity and integrals imply that
	$u_k,v_k$ are non-negative and non-increasing and that both have mass
	$m_k$.
	
	Now consider $I_k=1$.  The first part of
	\eqref{eq:coupling-condition}, together with the induction hypothesis,
	gives
	\[
	\left(\frac{u_{k-1}}{m_{k-1}}\right)(0+)
	\le\frac{h_{k-1}}{m_{k-1}}\le A,
	\]
	and the same estimate holds for $v_{k-1}/m_{k-1}$.  These normalized
	functions are non-negative and non-increasing and have integral one;
	hence
	\[
	\frac{u_{k-1}}{m_{k-1}},\quad
	\frac{v_{k-1}}{m_{k-1}}\in\Dcone_A.
	\]
	The second part of the condition says that $G_k$ occurs.  Applying
	Lemma~\ref{lem:lower} to both normalized functions and multiplying by
	$m_{k-1}$ gives the common lower bounds
	\begin{equation}\label{eq:common-lower-bound}
		\mathcal P_ku_{k-1}\ge\kappa_\ell m_{k-1},
		\qquad
		\mathcal P_kv_{k-1}\ge\kappa_\ell m_{k-1}.
	\end{equation}
	Therefore the subtractions in \eqref{eq:coupling-recursion} leave
	non-negative functions.  Subtracting a constant preserves their
	non-increasing property.  Since $m([0,1])=1$,
	\[
	\int u_k\,dm
	=m_{k-1}-\kappa_\ell m_{k-1}
	=(1-\kappa_\ell)m_{k-1}=m_k,
	\]
	and the same calculation applies to $v_k$.  Also
	$m_k=(1-\kappa_\ell)m_{k-1}>0$.
	
	It remains to check the height bound, in either case.  By
	\eqref{eq:height} and the induction hypothesis,
	\begin{align*}
		(\mathcal P_ku_{k-1})(0+)
		&\le\widehat S_k u_{k-1}(0+)
		+M\widehat R_k\int u_{k-1}\,dm\\
		&\le\widehat S_kh_{k-1}+M\widehat R_km_{k-1}=h_k.
	\end{align*}
	The same bound holds with $v_{k-1}$ in place of $u_{k-1}$.  If
	$I_k=0$, this is already the required estimate.  If $I_k=1$, the
	subtracted constant is non-negative and can only lower the value at
	$0+$.  Hence $u_k(0+),v_k(0+)\le h_k$ in both cases.  This proves
	\eqref{eq:coupling-invariants} at time $k$.  The base step and the
	induction step therefore define the coupling construction for every
	$k\ge1$.
	
	We next explain how the uncoupled mass controls the distance between the original evolved densities.
	On a successful coupling block, i.e., when $I_k=1$, define the normalized remainders
	\[
	\widetilde u_k=\frac{u_k}{m_k},\qquad
	\widetilde v_k=\frac{v_k}{m_k}.
	\]
	Because $m_k=(1-\kappa_\ell)m_{k-1}$ on this block, the recursion can be
	rewritten as
	\[
	\mathcal P_k\!\left(\frac{u_{k-1}}{m_{k-1}}\right)
	=\kappa_\ell+(1-\kappa_\ell)\widetilde u_k,
	\]
	\[
	\mathcal P_k\!\left(\frac{v_{k-1}}{m_{k-1}}\right)
	=\kappa_\ell+(1-\kappa_\ell)\widetilde v_k.
	\]
	This is the analogue of \cite[Equation~(13)]{StenlundSulku2014}:
	a fraction $\kappa_\ell$ of the two normalized densities has been made
	identical and removed from further comparison.
	
	Finally, subtracting the same common function from both evolved parts
	does not change their difference.  Inductively,
	\begin{align*}
		u_k-v_k
		&=\mathcal P_k(u_{k-1}-\kappa_\ell m_{k-1}I_{k})-\mathcal P_k(v_{k-1}-\kappa_\ell m_{k-1}I_{k})\\
		&=\mathcal P_k(u_{k-1}-u_{k-1}).
	\end{align*}
	Therefore, it follows that $u_k-v_k=\Lop_\omega^{k\ell}(f-g)$.
	Since $u_k,v_k\ge0$ and both have mass $m_k$,
	\begin{equation}\label{eq:mass-controls-distance}
		\norm{\Lop_\omega^{k\ell}(f-g)}_1
		=\norm{u_k-v_k}_1
		\le\norm{u_k}_1+\norm{v_k}_1
		=2m_k,
	\end{equation}
	where $\Vert f \Vert_1 = \int |f| \, dm$.
	Thus the uncoupled mass controls the distance between the two evolved
	densities.
	
	\begin{remark}
		Let
		$
		N_k=\sum_{j=1}^k I_j
		$
		be the number of successful couplings completed by the end of block
		$k$.  Each successful block multiplies the uncoupled mass by
		$1-\kappa_\ell$, while an unsuccessful block leaves it unchanged.  Hence
		\[
		m_k=(1-\kappa_\ell)^{N_k},
		\qquad
		\norm{\Lop_\omega^{k\ell}(f-g)}_1
		\le2(1-\kappa_\ell)^{N_k}.
		\]
		This is the fixed-block counterpart of the estimate in
		\cite{StenlundSulku2014} in terms of the number of couplings. We will not estimate the coupling times or the random variable $N_k$ directly.
		Instead, in Proposition~\ref{prop:conditional-contraction},
		we prove that a weighted sum $V_k=m_k+\varepsilon h_k$ contracts
		in conditional expectation from one block to the next.
	\end{remark}
	
	\subsection{Conditional block estimates}
	We now use
	independence of successive blocks to estimate, conditionally
	on the past, the probability of a successful block and the evolution
	of the height bound.  These estimates are
	needed in the contraction argument of Section~4.
	
	The coordinates in the $k$th block are independent and each has law
	$\eta$.  Hence
	\[
	\E\widehat S_k
	=\prod_{j=(k-1)\ell+1}^{k\ell}\E(\lambda_{\omega_j}^{-1})
	=\theta^\ell=q.
	\]
	Markov's inequality and the choice of $q$ therefore give
	\[
	\Pp(G_k^c)
	=\Pp\left(\widehat S_k>\frac{a}{2A}\right)
	\le\frac{2A}{a}\E\widehat S_k
	=\frac{2A}{a}q\le\frac12.
	\]
	Thus $\Pp(G_k)\ge1/2$.  We denote the current-block sigma-algebra by
	\[
	\mathcal G_k
	=\sigma(\omega_{(k-1)\ell+1},\ldots,\omega_{k\ell}).
	\]
	The variables $(\widehat S_k,\widehat R_k)$ and the event $G_k$ are
	$\mathcal G_k$-measurable, while $\mathcal G_k$ is independent of
	$\mathcal F_{k-1}$ under the product law.  Consequently,
	\begin{equation}\label{eq:conditional-good}
		\Pp(G_k\mid\mathcal F_{k-1})=\Pp(G_k)\ge\frac12
		\quad\text{a.s.}
	\end{equation}
	It also implies that the conditional expectations of the block
	variables equal their ordinary expectations.  In detail,
	\begin{align}
		\E(\widehat S_k\mid\mathcal F_{k-1})
		&=\E\widehat S_k=q,\label{eq:conditional-S}\\
		\E(\widehat R_k\mid\mathcal F_{k-1})
		&=\sum_{i=(k-1)\ell+1}^{k\ell}
		\E\left(\prod_{j=i+1}^{k\ell}\lambda_{\omega_j}^{-1}\right)\notag\\
		&=\sum_{i=(k-1)\ell+1}^{k\ell}\theta^{k\ell-i}
		=\sum_{r=0}^{\ell-1}\theta^r
		\le\frac1{1-\theta}.\label{eq:conditional-R}
	\end{align}
	
	\begin{lemma}\label{lem:conditional-height}
		Both $m_k$ and $h_k$  are
		$\mathcal F_k$-measurable and satisfy
		\[
		0<m_k\le1,\qquad h_k\ge0,\qquad \E h_k<\infty.
		\]
		Moreover, for every $k\ge1$,
		\begin{equation}\label{eq:conditional-height}
			\E(h_k\mid\mathcal F_{k-1})
			\le qh_{k-1}+Dm_{k-1}.
		\end{equation}
	\end{lemma}
	\begin{proof}
		The measurability and pointwise bounds follow by induction from
		\eqref{eq:first-coupling-recursion} and
		\eqref{eq:coupling-recursion}; here $\kappa_\ell\le1/2$ ensures that
		$m_k$ remains strictly positive.  Integrability follows from the same
		induction and the conditional expectations in
		\eqref{eq:conditional-S}-\eqref{eq:conditional-R}.
		Finally, using the $\mathcal F_{k-1}$-measurability of
		$m_{k-1}$ and $h_{k-1}$ together with
		\eqref{eq:conditional-S} and \eqref{eq:conditional-R}, we obtain
		\begin{align*}
			\E(h_k\mid\mathcal F_{k-1})
			&=\E\bigl(\widehat S_kh_{k-1}
			+M\widehat R_km_{k-1}\mid\mathcal F_{k-1}\bigr)\\
			&=h_{k-1}\E(\widehat S_k\mid\mathcal F_{k-1})
			+Mm_{k-1}\E(\widehat R_k\mid\mathcal F_{k-1})\\
			&\le qh_{k-1}+Dm_{k-1},
		\end{align*}
		as claimed.
	\end{proof}
	
	\section{Conditional contraction and exponential memory loss}
	The purpose of this section is to turn the coupling construction of the preceding section into a uniform contraction estimate and thereby prove exponential loss of memory. The uncoupled mass $m_k$ directly controls the distance between the two evolved densities.
	But $m_k$ need not decrease on every block: when $h_{k-1}>Am_{k-1}$, no coupling is performed and the mass remains unchanged. The height $h_k$ contracts in conditional expectation under the local average expansion assumption, its estimate contains an additional term proportional to $m_{k-1}$ by Lemma \ref{lem:conditional-height}. Thus neither $m_k$ nor $h_k$ alone admits a uniform one-block contraction. So we introduce
	$
	V_k:=m_k+\varepsilon h_k,
	$
	where $\varepsilon>0$ is chosen sufficiently small. In the low-height regime $h_{k-1}\le Am_{k-1}$, a coupling block occurring with uniformly positive conditional probability removes a fixed fraction of $m_{k-1}$, and this decrease dominates any possible increase of the height term. In the high-height regime $h_{k-1} > Am_{k-1}$, coupling is temporarily unavailable, but the conditional height estimate
	\eqref{eq:conditional-height}
	supplies the required contraction. The main objective of this section is therefore to prove the following  conditional contraction estimate. Since $m_k\leq V_k$, this conditional contraction will then yield both annealed and quenched exponential loss of memory.

	\subsection{One-block conditional contraction}
	
	\begin{proposition}
		\label{prop:conditional-contraction}
		Let
		\[
		V_k=m_k+\varepsilon h_k,
		\qquad
		\varepsilon=\frac{\kappa_\ell}{4(D+1)}.
		\]
		Then there exists a constant $r\in(0,1)$, depending
		only on the constants fixed in Section~3, such that
		\begin{equation}\label{eq:conditional-contraction}
			\E(V_k\mid\mathcal F_{k-1})\le rV_{k-1},
			\qquad k\ge1.
		\end{equation}
		Consequently,
		\begin{equation}\label{eq:mean-contraction}
			\E V_k\le (1+\varepsilon B)r^k,
			\qquad k\ge0.
		\end{equation}
	\end{proposition}
	\begin{proof}
		We divided the following two case:
		
		\emph{Low-height regime.}
		First suppose that $h_{k-1}\le Am_{k-1}$.  This condition is
		$\mathcal F_{k-1}$-measurable, and conditional on the past a coupling
		occurs precisely when the independent event $G_k$ occurs.  Since
		$\Pp(G_k\mid\mathcal F_{k-1})\ge1/2$ by
		\eqref{eq:conditional-good}, we have $I_k=\one_{G_k}$ in this regime.
		As $m_{k-1}$ is $\mathcal F_{k-1}$-measurable, the mass recursion gives
		\begin{align}
			\E(m_k\mid\mathcal F_{k-1})
			&=\E\bigl(m_{k-1}(1-\kappa_\ell I_k)
			\mid\mathcal F_{k-1}\bigr)\notag\\
			&=m_{k-1}\left(1-\kappa_\ell
			\E(I_k\mid\mathcal F_{k-1})\right)\notag\\
			&=m_{k-1}\left(1-\kappa_\ell
			\Pp(G_k\mid\mathcal F_{k-1})\right)\notag\\
			&\le
			\left(1-\frac{\kappa_\ell}{2}\right)m_{k-1}.
			\label{eq:conditional-mass-low}
		\end{align}
		Thus the low-height regime produces a uniform expected contracting of the
		uncoupled mass.  Combining this with
		\eqref{eq:conditional-height},
		\begin{align*}
			\E(V_k\mid\mathcal F_{k-1})
			&=\E(m_k\mid\mathcal F_{k-1})
			+\varepsilon\E(h_k\mid\mathcal F_{k-1})\\
			&\le
			\left(1-\frac{\kappa_\ell}{2}\right)m_{k-1}
			+\varepsilon(qh_{k-1}+Dm_{k-1})\\
			&=V_{k-1}-\left(\frac{\kappa_\ell}{2}-\varepsilon D\right)m_{k-1}-\varepsilon(1-q)h_{k-1}\\
			&\le
			V_{k-1}-\frac{\kappa_\ell}{4}m_{k-1}
			\le r_1V_{k-1},
		\end{align*}
		where
		$
		r_1=1-\frac{\kappa_\ell}{4(1+\varepsilon A)}<1.
		$
		Indeed, $\varepsilon D\le\kappa_\ell/4$, while the present regime gives
		$V_{k-1}\le(1+\varepsilon A)m_{k-1}$.  The last two lines above follow from these
		inequalities.
		
		\emph{High-height regime.}
		Now suppose that $h_{k-1}>Am_{k-1}$.  Then $I_k=0$ and
		$m_k=m_{k-1}$, so there is no mass contraction in this block.  However,
		$m_{k-1}<h_{k-1}/A$.  Since $q\le1/4$ and $D/A\le1/8$,
		\[
		(1-q)h_{k-1}-Dm_{k-1}
		\ge\left(1-q-\frac DA\right)h_{k-1}
		\ge\frac12h_{k-1}.
		\]
		It follows that
		\begin{align*}
			\E(V_k\mid\mathcal F_{k-1})
			&=m_{k-1}+\varepsilon
			\E(h_k\mid\mathcal F_{k-1})\\
			&\le m_{k-1}
			+\varepsilon(qh_{k-1}+Dm_{k-1})\\
			&=
			V_{k-1}-\varepsilon\bigl((1-q)h_{k-1}-Dm_{k-1}\bigr)\\
			&\le
			V_{k-1}-\frac{\varepsilon}{2}h_{k-1} 
			\le r_2V_{k-1},
		\end{align*}
		where
		$
		r_2=1-\frac{\varepsilon}{2(\varepsilon+1/A)}<1.
		$
		Here the last inequality follows from
		$V_{k-1}\le(\varepsilon+1/A)h_{k-1}$.
		
		Thus, with
		$r=\max\{r_1,r_2\}<1$,
		$
		\E(V_k\mid\mathcal F_{k-1})\le rV_{k-1}.
		$
		Because $V_{k-1}$ is $\mathcal F_{k-1}$-measurable and integrable, it follows that
		\[
		\E V_k
		=\E\!\left[\E(V_k\mid\mathcal F_{k-1})\right]
		\le r\E V_{k-1}.
		\]
		Iteration of this ordinary expectation inequality, together with
		$V_0=1+\varepsilon B$, gives
		\[
		\E V_k\le r^k(1+\varepsilon B).
		\]
		This proves both \eqref{eq:conditional-contraction} and
		\eqref{eq:mean-contraction}.
	\end{proof}
	
	\subsection{Annealed and quenched exponential memory loss}
	\begin{proof}[Completion of the proof of Theorem~\ref{thm:main}]
		By $m_k\le V_k$, \eqref{eq:mass-controls-distance}, and \eqref{eq:mean-contraction}
		it follows that
		\[
		\E
		\norm{\Lop_\omega^{k\ell}(f-g)}_1
		\le2(1+\varepsilon B)r^k.
		\]
		For $n=k\ell+j$, $0\le j<\ell$, the $L^1$ contraction property gives
		\[
		\E
		\norm{\Lop_\omega^n(f-g)}_1
		\le2(1+\varepsilon B)r^{\lfloor n/\ell\rfloor}.
		\]
		Since $\lfloor n/\ell\rfloor\ge n/\ell-1$, the choices
		\[
		\rho=r^{1/\ell},
		\qquad
		C_B=2r^{-1}(1+\varepsilon B)
		\]
		give
		\[
		\E\norm{\Lop_\omega^n(f-g)}_1\le C_B\rho^n.
		\]
		This is the annealed estimate in Theorem~\ref{thm:main}.
		
		Choose
		$
		0<c<\frac1\ell\log\frac1r,
		$
		so that $e^{c\ell}r<1$, and set
		\[
		C_{B,c}(\omega)
		=2e^{c\ell}\sum_{k=0}^\infty e^{ck\ell}V_k(\omega).
		\]
		All summands are non-negative.  The estimate for
		$\E V_k$ therefore yield
		\begin{align*}
			\E C_{B,c}
			&=2e^{c\ell}\E\left(\sum_{k=0}^\infty
			e^{ck\ell}V_k\right)
			=2e^{c\ell}\sum_{k=0}^\infty e^{ck\ell}\E V_k
			&\le2e^{c\ell}(1+\varepsilon B)
			\sum_{k=0}^\infty(e^{c\ell}r)^k
			&=\frac{2e^{c\ell}(1+\varepsilon B)}{1-e^{c\ell}r}<\infty.
		\end{align*}
		Thus $C_{B,c}\in L^1(\Pp)$; in particular,
		$C_{B,c}(\omega)<\infty$ for almost every $\omega$.
		
		Fix such  $\omega$ and write $n=k\ell+j$, where $0\le j<\ell$.
		The $L^1$ contraction of the remaining $j$ transfer operators and
		\eqref{eq:mass-controls-distance} give
		\[
		\norm{\Lop_\omega^n(f-g)}_1
		\le\norm{\Lop_\omega^{k\ell}(f-g)}_1
		\le2m_k\le2V_k.
		\]
		Since $e^{cj}\le e^{c\ell}$,
		\begin{align*}
			e^{cn}\norm{\Lop_\omega^n(f-g)}_1
			&\le2e^{c\ell}e^{ck\ell}V_k \le2e^{c\ell}\sum_{p=0}^\infty e^{cp\ell}V_p
			=C_{B,c}(\omega).
		\end{align*}
		Consequently,
		\[
		\norm{\Lop_\omega^n(f-g)}_1
		\le C_{B,c}(\omega)e^{-cn},
		\qquad n\ge0.
		\]
		
	\end{proof}

	\section{Consequences for bounded-variation densities and the averaged operator}
	
	We first extend Theorem~\ref{thm:main} from decreasing probability
	densities to arbitrary real $BV$ functions of zero integral.  We then
	use this extension to identify the unique absolutely continuous
	stationary distribution of the annealed dynamics.
	
	\begin{lemma}
		\label{lem:bvdecomp}
		Let $\psi\in BV([0,1])$ be real-valued, non-zero in $L^1$, and satisfy
		$\int \psi\,dm=0$.
		Then there exist $f_{\psi},g_{\psi}\in\Dcone_3$ and a constant
		$c_\psi\in(0,2\|\psi\|_{\mathrm{BV}}]$ such that
		\[
		\psi=c_\psi(f_{\psi}-g_{\psi})\qquad\text{a.e.}
		\]
	\end{lemma}
	
	\begin{proof}
		Choose the representative $\bar\psi$ for $\psi \in BV([0,1])$
		that is right-continuous
		on $(0,1)$, with $\bar\psi(0)= \bar\psi(0+)$ and
		$\bar\psi(1)= \bar\psi(1-)$. We write $V_{[a,b]}(\varphi)$
		for the usual variation of a function $\varphi : [a,b] \to \mathbb{R}$:
		\begin{align*}
			V_{[a,b]}(\varphi) =
			\sup \biggl\{
			\sum_{j=0}^{t-1} | \varphi( x_{j+1} ) -  \varphi( x_{j} )  | \: : \:
			a = x_0 < \ldots < x_t = b
			\biggr\}
		\end{align*}
		We have
		\[
		V_{  [0,1] }(\bar\psi)=\operatorname{Var}(\psi),
		\qquad
		\sup_{x\in[0,1]}|\bar\psi(x)|=\|\psi\|_{L^\infty}.
		\]

		Let
		\[
		W(x)=V_{[x,1]}(\bar\psi),
		\qquad
		f=\|\psi\|_{\mathrm{BV}}+W,
		\qquad
		g=\|\psi\|_{\mathrm{BV}}+W-\bar\psi.
		\]
		The function $W$ is non-negative and non-increasing, and hence so is
		$f$.  Moreover, if $x<y$, additivity of variation gives
		\[
		g(x)-g(y)
		=V_{[x,y]}(\bar\psi)-\bar\psi(x)+\bar\psi(y)\ge0,
		\]
		so $g$ is also non-increasing.  Both functions are non-negative since
		$\|\psi\|_{\mathrm{BV}}\ge\norm{\psi}_{L^\infty}$.
		
		Because $f-g=\bar\psi$ and $\int \bar\psi\,dm=0$, the functions $f$ and $g$ have the
		same integral. Let
		\[
		c_\psi:= \int_0^1 f \,dm= \int_0^1 (\|\psi\|_{\mathrm{BV}}+W)\,dm.
		\]
		Since $0\le W\le
		\operatorname{Var}(\psi)$, it follows that
		$
		\|\psi\|_{\mathrm{BV}}\le c_\psi
		\le \|\psi\|_{\mathrm{BV}}+\operatorname{Var}(\psi)\le2\|\psi\|_{\mathrm{BV}}.
		$
		Furthermore,
		\[
		f(0+)\le2\|\psi\|_{\mathrm{BV}},
		\qquad
		g(0+)\le \|\psi\|_{\mathrm{BV}}+\operatorname{Var}(\psi)+\norm{\psi}_{L^\infty}=2\|\psi\|_{\mathrm{BV}}.
		\]
		Consequently, the normalized functions
		\[
		f_\psi=\frac{f}{c_\psi},
		\qquad
		g_\psi=\frac{g}{c_\psi}
		\]
		are non-negative, non-increasing probability densities with height at
		most $2$, and hence belong to $\Dcone_3$.  Finally,
		$$
		\psi
		= \bar{\psi} =
		f-g=c_\psi(f_\psi-g_\psi) \quad \text{a.e.}
		$$
	\end{proof}
	
	\begin{corollary}
		\label{cor:bv-memory}
		There exist constants $C<\infty$, $\rho\in(0,1)$, and
		$c>0$, and a non-negative random variable $Z\in L^1(\Pp)$, such that
		for every $\psi\in BV([0,1])$ with $\int \psi\,dm=0$,
		\begin{equation}\label{eq:bv-annealed}
			\E\norm{\Lop_\omega^n\psi}_1
			\le C\|\psi\|_{\mathrm{BV}}\rho^n,
			\qquad n\ge0,
		\end{equation}
		and, almost surely
		\begin{equation}\label{eq:bv-quenched}
			\norm{\Lop_\omega^n\psi}_1
			\le Z(\omega)\|\psi\|_{\mathrm{BV}} e^{-cn},
			\qquad n\ge0.
		\end{equation}
	\end{corollary}
	
	\begin{proof}
		The assertion is immediate for $\psi=0$.  Otherwise, apply
		Lemma~\ref{lem:bvdecomp} and write
		\[
		\psi=c_\psi(f_\psi-g_\psi),
		\qquad
		f_\psi,g_\psi\in\Dcone_3,
		\qquad
		c_\psi\le2\|\psi\|_{\mathrm{BV}}.
		\]
		The annealed conclusion of Theorem~\ref{thm:main}, with $B=3$, yields
		\[
		\E\norm{\Lop_\omega^n\psi}_1
		\le c_\psi C_3\rho^n
		\le2C_3\|\psi\|_{\mathrm{BV}}\rho^n.
		\]
		The quenched conclusion of the same theorem gives, on the common
		full-measure event associated with $B=3$,
		\[
		\norm{\Lop_\omega^n\psi}_1
		\le c_\psi C_{3,c}(\omega)e^{-cn}
		\le2C_{3,c}(\omega)\|\psi\|_{\mathrm{BV}} e^{-cn}.
		\]
		Thus one may take $C=2C_3$ and $Z=2C_{3,c}$.
	\end{proof}
	
	Define the averaged transfer operator by
	\[
	\Qop f=\E(L_{\omega_1}f).
	\]
	It is a positive, integral-preserving contraction on $L^1(m)$.
	
	\begin{proposition}
		\label{prop:annealed-density}
		The operator $\Qop$ has a unique invariant probability density $\phi$ in
		$L^1(m)$.   Then
		$\phi\in BV\cap L^\infty$ and $\phi$ is non-increasing.  Moreover, for
		every $BV$ probability density $f$,
		\begin{equation}\label{eq:annealed-bv-convergence}
			\norm{\Qop^nf-\phi}_1
			\le C  \|f-\phi\|_{\mathrm{BV}}\rho^n,
			\qquad n\ge0,
		\end{equation}
		where $C$ and $\rho$ are as in Corollary~\ref{cor:bv-memory}.  For
		every probability density $f\in L^1(m)$,
		\begin{equation}\label{eq:annealed-l1-convergence}
			\norm{\Qop^nf-\phi}_1\longrightarrow0,
			\qquad n\to\infty.
		\end{equation}
	\end{proposition}
	
	\begin{proof}
		Independence of the  coordinates  of $\omega$ and Fubini's theorem give
		\begin{equation}\label{eq:averaged-iterate}
			\Qop^nf=\E(\Lop_\omega^nf),
			\qquad f\in L^1(m),\quad n\ge1.
		\end{equation}
		In particular, positivity of the random transfer operators implies
		that $\Qop$ is an $L^1$ contraction.
		
		\emph{Step 1: an invariant cone.}
		Recall $D=\frac{M}{1-\theta}$, and  set
		$
		B_0=\max\{1,D\}.
		$
		The operator $\Qop$ preserves non-negativity, integrals, and
		monotonicity.  If $f\in\Dcone_{B_0}$, monotone convergence at $0+$ and
		Lemma~\ref{lem:onestep} yield
		\[
		(\Qop f)(0+)
		\le\theta f(0+)+M
		\le\theta B_0+M
		\le B_0,
		\]
		where the last inequality follows from
		$M=(1-\theta)D\le(1-\theta)B_0$.  Hence
		\begin{equation}\label{eq:averaged-invariant-cone}
			\Qop(\Dcone_{B_0})\subset\Dcone_{B_0}.
		\end{equation}
		
		\emph{Step 2: construction of an invariant density.}
		Fix $f_0\in\Dcone_{B_0}$.  By
		\eqref{eq:averaged-invariant-cone}, both $f_0$ and $\Qop^jf_0$ belong
		to $\Dcone_{B_0}$ for every $j\ge0$.  Using
		\eqref{eq:averaged-iterate}, Jensen's inequality, and
		Theorem~\ref{thm:main}, we obtain
		\begin{align*}
			\norm{\Qop^{n+j}f_0-\Qop^nf_0}_1
			=\norm{\Qop^n(\Qop^jf_0-f_0)}_1 
			\le\E\norm{\Lop_\omega^n(\Qop^jf_0-f_0)}_1 
			\le C_{B_0}\rho^n.
		\end{align*}
		The bound is uniform in $j$.  Thus $(\Qop^nf_0)_{n\ge0}$ is a Cauchy
		sequence in $L^1(m)$; let
		\[
		\phi=\lim_{n\to\infty}\Qop^nf_0.
		\]
		
		The cone $\Dcone_{B_0}$ is closed in $L^1$.  Indeed, from any
		$L^1$-convergent sequence in this cone one may extract an almost
		everywhere convergent subsequence, and its limit admits a
		non-increasing representative with integral one and height at most
		$B_0$.  Therefore $\phi\in\Dcone_{B_0}$.  Since $\Qop$ is continuous
		on $L^1$,
		\[
		\Qop\phi
		=\lim_{n\to\infty}\Qop^{n+1}f_0
		=\phi.
		\]
		Thus $\phi$ is an invariant probability density. Since $\phi\in \Dcone_{B_0}$,  it follows that  it is bounded and non-increasing; in
		particular, $\phi\in BV\cap L^\infty$.
		
		\emph{Step 3: exponential attraction of $BV$ densities.}
		Let $f$ be a $BV$ probability density.  Then $f-\phi\in BV$ and
		$\int(f-\phi)\,dm=0$.  By \eqref{eq:averaged-iterate} and
		Corollary~\ref{cor:bv-memory},
		\begin{align*}
			\norm{\Qop^nf-\phi}_1
			=\norm{\Qop^n(f-\phi)}_1 
			\le\E\norm{\Lop_\omega^n(f-\phi)}_1 
			\le C \|f-\phi\|_{\mathrm{BV}}\rho^n,
		\end{align*}
		which proves \eqref{eq:annealed-bv-convergence}.
		
		\emph{Step 4: attraction of arbitrary $L^1$ densities.}
		Let $f\in L^1(m)$ be a probability density.  Choose non-negative
		step-function probability densities $(f_j)_{j\ge1}$ such that
		$\norm{f-f_j}_1\to0$.  Each $f_j$ belongs to $BV$, and the
		$L^1$-contractivity of $\Qop$, together with Step~3, gives
		\begin{align*}
			\limsup_{n\to\infty}\norm{\Qop^nf-\phi}_1
			&\le\norm{f-f_j}_1
			+\limsup_{n\to\infty}\norm{\Qop^nf_j-\phi}_1\\
			&=\norm{f-f_j}_1.
		\end{align*}
		Letting $j\to\infty$ proves \eqref{eq:annealed-l1-convergence}.
		
		\emph{Step 5: uniqueness.}
		If $\psi$ is an invariant probability density in $L^1$, then
		$\Qop^n\psi=\psi$ for all $n$.  Applying
		\eqref{eq:annealed-l1-convergence} with $f=\psi$ gives
		$\norm{\psi-\phi}_1=0$.  Hence $\phi$ is the unique invariant
		probability density in $L^1(m)$.
	\end{proof}
	
	\begin{remark}
		\label{rem:annealed-correlations}
		Let $f\in BV$ be a probability density and let $g\in L^\infty$.  By
		duality and \eqref{eq:annealed-bv-convergence},
		\[
		\left|
		\E\int g\circ T_\omega^n\,f\,dm
		-\int g\phi\,dm
		\right|
		\le C\|f-\phi\|_{\mathrm{BV}}\norm{g}_\infty\rho^n.
		\]
		This estimate means that the
		time-$n$ expectation issued from $f$ converges  to  the stationary expectation at an exponential rate.

		For correlations with respect to the stationary density $\phi$, let $f\in BV$ be real-valued,
		let $g\in L^\infty$, and put
		\[
		\bar f=\int f\phi\,dm,
		\qquad
		\psi=(f-\bar f)\phi.
		\]
		Since $\phi\in BV\cap L^\infty$, one has $\psi\in BV$, and
		$\int \psi\,dm=0$.  Duality and Corollary~\ref{cor:bv-memory} therefore
		give
		\begin{align*}
			&\left|
			\E\int f\,(g\circ T_\omega^n)\phi\,dm
			-\int f\phi\,dm\int g\phi\,dm
			\right|\\
			&\qquad
			=\left|\E\int g\,\Lop_\omega^n\psi\,dm\right|\\
			&\qquad
			\le C\|\psi\|_{\mathrm{BV}}\norm{g}_\infty\rho^n.
		\end{align*}
		This is exponential decay of annealed correlations for the stationary
		density $\phi$.
	\end{remark}

	\section{Applications to limit theorems}
	\label{sec:quenched-asip}
	
	We return to the i.i.d.\ setting of
	Theorem~\ref{thm:main}. We show that exponential memory loss,
	together with the preservation of decreasing functions, allows us
	to apply the abstract theorems in \cite{DHL2026},
	yielding a quenched ASIP and a rate of convergence in the quenched CLT
	for random Lasota--Yorke convex maps
	satisfying the local average expansion condition.
	To this end, we construct a family of equivariant
	sample measures and establish an exponential $L^s$ memory-loss
	estimate for products of a Lipschitz observable and
	a Lipschitz observable evolved by the normalized transfer operators.
	The random prefactor in this estimate has finite moments of
	every order.

	We work on the two-sided product extension
	\[
	(\widehat\Omega,\mathbb P,\sigma)
	=
	(\Omega^{\mathbb Z},\eta^{\otimes\mathbb Z},\sigma),
	\]
	where $\sigma$ is a full shift. For $\omega\in\widehat\Omega$, write $T_\omega=T_{\omega_0}$ and
	let $ L_\omega$ be its transfer operator with respect to
	Lebesgue measure $m$. Set
	\[
	T_\omega^n
	=
	T_{\sigma^{n-1}\omega}\circ\cdots\circ T_\omega,
	\qquad
	\mathcal L_\omega^n
	=
	L_{\sigma^{n-1}\omega}\cdots L_\omega.
	\]

	We recall two consequences of the preceding sections. First,
	$L_\omega$ preserves non-negativity, integrals, and
	monotonicity, and
	\begin{equation}
		\label{eq:asip-one-step}
		(L_\omega f)(0+)
		\le q(\omega)f(0+)+M\int f\,dm,
		\qquad
		q(\omega)=\lambda_{\omega_0}^{-1},
		\qquad
		\mathbb E q=\theta<1,
	\end{equation}
	for every non-negative non-increasing function $f$.
	Second, Corollary~\ref{cor:bv-memory} gives constants
	$C<\infty$, $\rho\in(0,1)$, $c>0$, and a non-negative
	$Z\in L^1(\mathbb P)$ such that, whenever $\int \psi\,dm=0$,
	\begin{equation}
		\label{eq:asip-memory-input}
		\mathbb E\|\mathcal L_\omega^n \psi\|_1
		\le C\rho^n \|\psi\|_{\mathrm{BV}},
		\qquad
		\|\mathcal L_\omega^n \psi\|_1
		\le Z(\omega)e^{-cn}\|\psi\|_{\mathrm{BV}}.
	\end{equation}
	The second estimate holds on a common full measure set, for all such $\psi$ and all $n\geq 0$.
	By its construction, $Z(\omega)$ is measurable with respect to $\sigma(  \omega_j \: : \: j \ge 0 )$.
	
	\begin{lemma}
		\label{lem:asip-sample-densities}
		There exists a measurable subset $\Omega_* \subset \Omega$ with $\mathbb{P}(\Omega_*) = 1$ and
		a
		measurable map $h : \Omega \to L^1(m)$ such that the following hold for every $\omega \in \Omega_*$:
		\begin{itemize}
			\item $h_\omega = h(\omega) \in BV([0,1])$ is non-increasing and $\int h_\omega \, dm = 1$;
			\item The equivariance property holds:
			$$
			L_\omega h_\omega=h_{\sigma\omega};
			$$
			\item The map $h$ is $\sigma( \omega_j \: : \:  j < 0 )$ measurable and
			\begin{equation}
				\label{eq:asip-sample-height}
				h_\omega(0+)\le H(\omega),
				\qquad
				H(\omega)
				=
				1+M\sum_{j=0}^{\infty}
				\prod_{t=1}^{j}q(\sigma^{-t}\omega),
				\qquad
				\mathbb EH=1+\frac{M}{1-\theta}.
			\end{equation}
		\end{itemize}
	\end{lemma}
	
	\begin{proof}
		Define the pullback densities
		\[
		h_\omega^{(0)}=\mathbf1,
		\qquad
		h_\omega^{(n)}
		=
		\mathcal L_{\sigma^{-n}\omega}^{n}\mathbf1,
		\qquad n\ge1.
		\]
		Each $h_\omega^{(n)}$ is a non-negative decreasing probability
		density.
		
		We first prove the  convergence of $h_\omega^{(n)}$  in $L^1(m)$. Notice that
		$\mathcal L_\omega\mathbf1$ is non-negative and
		non-increasing, and \eqref{eq:asip-one-step} gives
		$L_\omega\mathbf1(0+)\le q(\omega)+M$. Hence
		\[
		\|L_\omega\mathbf1-\mathbf1\|_\infty
		\le1+q(\omega)+M,
		\qquad
		\operatorname{Var}(L_\omega\mathbf1-\mathbf1)
		=
		\operatorname{Var}(L_\omega\mathbf1)
		\le q(\omega)+M.
		\]
		Consequently,
		\begin{equation}
			\label{eq:asip-first-increment}
			\| L_\omega\mathbf1-\mathbf1\|_{\mathrm{BV}}
			\le1+2q(\omega)+2M,
			\qquad
			\mathbb E\mathfrak \| L_\omega\mathbf1-\mathbf1\|_{\mathrm{BV}}
			\le1+2\theta+2M<\infty.
		\end{equation}
		Notice that the difference of two successive pullback densities is
		\[
		h_\omega^{(n+1)}-h_\omega^{(n)}
		=
		\mathcal L_{\sigma^{-n}\omega}^{n}
		\bigl( L_{\sigma^{-n-1}\omega}\mathbf1-\mathbf1\bigr).
		\]
		The function $ L_{\sigma^{-n-1}\omega}\mathbf1-\mathbf1 $  has zero integral and depends only
		on $\omega_{-n-1}$, whereas the subsequent $n$ operators depend
		on $\omega_{-n},\ldots,\omega_{-1}$. These coordinates are
		independent. Using the second estimate in \eqref{eq:asip-memory-input}, independence, and
		\eqref{eq:asip-first-increment}, we obtain
		\begin{align*}
			\mathbb E\|h_\omega^{(n+1)}-h_\omega^{(n)}\|_1
			&\le \mathbb{E} [   Z( \sigma^{-n} \omega) e^{-cn} \|  L_{\sigma^{-n-1}\omega}\mathbf1-\mathbf1  \|_{\mathrm{BV}}  ] \\
			&= \mathbb{E} [   Z( \sigma^{-n}  \omega) ] \mathbb{E} [ \|  L_{\sigma^{-n-1}\omega}\mathbf1-\mathbf1  \|_{\mathrm{BV}}  ]
			e^{-cn}
			\le C 	e^{-cn}.
		\end{align*}
		Therefore
		\[
		\mathbb E\sum_{n=0}^{\infty}
		\|h_\omega^{(n+1)}-h_\omega^{(n)}\|_1<\infty.
		\]
		Set
		$$
		\Omega_* = \biggl\{   \omega \in \Omega \: : \:   \sum_{n=0}^{\infty}
		\|h_\omega^{(n+1)}-h_\omega^{(n)}\|_1<\infty  \biggr\}.
		$$
		Then $\mathbb{P}(\Omega_*) = 1$ and
		$( h_\omega^{(n)} )_{n \ge 0}$ is a Cauchy sequence in $L^1(m)$ for each $\omega \in \Omega_*$. We define
		$$
		h(\omega) = h_\omega = \lim_{n\to\infty} h_\omega^{(n)}  \in L^1(m), \quad \omega \in \Omega_*.
		$$
		Clearly $h_\omega$ is non-negative and
		\begin{align*}
			L_\omega h_\omega = L_\omega(  \lim_{n \to \infty}   h_\omega^{(n-1)}    )
			= \lim_{n \to \infty}  L_\omega L_{ \sigma^{-1} \omega   } \cdots L_{  \sigma^{ - ( n - 1) } \omega  } \mathbf{1}
			= \lim_{ n \to \infty }  \mathcal{L}_{  \sigma^{-n} \sigma \omega  }^n \mathbf{1} =
			h_{\sigma \omega}.
		\end{align*}
		Moreover, $h_\omega$ is non-increasing because the family of non-increasing elements in $L^1(m)$ is closed.
		
		From the joint measurability of $(\xi, x) \mapsto T_\xi(x)$ it follows that $\omega \mapsto h_\omega^{(n)}$ is a measurable map
		into $L^1(m)$. Therefore $\Omega_*$ is measurable. Moreover (see e.g. \cite[Lemma 1.10]{Kallenberg2002}), since
		$$
		\lim_{n \to \infty} \Vert  h_\omega^{(n)} - h_\omega  \Vert_1 = 0 \quad \text{for $\omega \in \Omega_*$},
		$$
		also $h : \Omega_* \to L^1(m)$ is measurable.

		As $h_\omega$ is non-increasing, we have
		$$
		\text{Var}(h_\omega) \le h_\omega(0+).
		$$
		It therefore remains to control $h_\omega(0+)$. To this end, observe that
		\[
		h_\omega^{(n)}
		=
		\mathcal L_{\sigma^{-1}\omega}
		h_{\sigma^{-1}\omega}^{(n-1)}.
		\]
		By \eqref{eq:asip-one-step}, it follows that
		\[
		h_\omega^{(n)}(0+)
		\le
		q(\sigma^{-1}\omega)
		h_{\sigma^{-1}\omega}^{(n-1)}(0+)+M.
		\]
		Starting with $h_\omega^{(0)}(0+)=1$ and iterating this
		inequality gives
		\begin{equation}
			\label{eq:asip-pullback-height}
			h_\omega^{(n)}(0+)
			\le
			\prod_{t=1}^{n}q(\sigma^{-t}\omega)
			+
			M\sum_{j=0}^{n-1}
			\prod_{t=1}^{j}q(\sigma^{-t}\omega).
		\end{equation}
		Independence gives
		\[
		\mathbb E\prod_{t=1}^{j}q(\sigma^{-t}\omega)=\theta^j.
		\]
		Consequently
		\[
		\mathbb E\sum_{j=0}^{\infty}
		\prod_{t=1}^{j} q(\sigma^{-t}\omega)
		=
		\sum_{j=0}^{\infty}\theta^j
		=
		\frac1{1-\theta}<\infty.
		\]
		In particular,
		$$
		\lim_{n \to \infty} \prod_{t=1}^{n}q(\sigma^{-t}\omega) = 0 \quad \text{for $\mathbb{P}$-a.e. $\omega \in \Omega$.}
		$$
		Passing to an almost everywhere convergent subsequence of
		$h_\omega^{(n)}$ in \eqref{eq:asip-pullback-height} shows that
		$h_\omega(0+)\le H(\omega)$, with $H$ as in
		\eqref{eq:asip-sample-height}.
	\end{proof}
	
	Set $d\mu_\omega=h_\omega\,dm$. The preceding lemma gives
	\[
	(T_\omega)_*\mu_\omega=\mu_{\sigma\omega}.
	\]
	Let $P_\omega$ denote the transfer operator relative to these
	sample measures. For bounded measurable $f$, it is given by
	\begin{equation}
		\label{eq:asip-normalized-operator}
		P_\omega f
		=
		\frac{ L_\omega(fh_\omega)}{h_{\sigma\omega}}
		\qquad\mu_{\sigma\omega}\text{-almost everywhere}.
	\end{equation}
	We define the quotient to be zero where the denominator
	vanishes. Note that
	\[
	|L_\omega(fh_\omega)|
	\le
	\|f\|_\infty L_\omega h_\omega
	=
	\|f\|_\infty h_{\sigma\omega},
	\]
	so the numerator also vanishes on that set.
	For every bounded measurable $g$, we have
	\[
	\int (P_\omega f)g\,d\mu_{\sigma\omega}
	=
	\int f(g\circ T_\omega)\,d\mu_\omega.
	\]
	Moreover,
	\[
	P_\omega\mathbf1=\mathbf1
	\quad\mu_{\sigma\omega}\text{-almost everywhere},
	\qquad
	\|P_\omega f\|_\infty\le\|f\|_\infty.
	\]

	Writing
	$P_\omega^n=P_{\sigma^{n-1}\omega}\cdots P_\omega$,
	we obtain by induction
	\begin{equation}
		\label{eq:asip-normalized-iterate}
		(P_\omega^nf)h_{\sigma^n\omega}
		=
		\mathcal L_\omega^n(fh_\omega).
	\end{equation}
	Moreover,
	\begin{equation}
		\label{eq:asip-normalized-l1}
		\|P_\omega^nf\|_{L^1(\mu_{\sigma^n\omega})}
		=
		\|\mathcal L_\omega^n(fh_\omega)\|_{L^1(m)}.
	\end{equation}
	For an integrable function $f$, we write
	$[f]_\omega=f-\mu_\omega(f)$.
	
	\begin{lemma}
		\label{lem:asip-two-observable}
		For every $1\le s,p<\infty$, there exist a deterministic
		constant $\gamma_{s,p}>0$, a non-negative random variable
		$K_{s,p}\in L^p(\mathbb P)$ and $\Omega_{s,p}^{*} \subset \Omega$ with $\Pp(\Omega_{s,p}^{*})=1$ such that the following hold: for almost every $\omega \in \Omega_{s,p}^{*}$,
		every pair of Lipschitz functions $g_1,g_2$, and all integers
		$r\le i\le j$,
		\begin{align}
			\label{eq:asip-two-observable}
			&\left\|
			P_{\sigma^i\omega}^{j-i}
			\left[g_2P_{\sigma^r\omega}^{i-r}g_1\right]_{\sigma^i\omega}
			\right\|_{L^s(\mu_{\sigma^j\omega})}
			\le
			K_{s,p}(\sigma^i\omega)e^{-\gamma_{s,p}(j-i)}
			(1+\|g_1\|_{\mathrm{Lip}})
			(1+\|g_2\|_{\mathrm{Lip}}).
		\end{align}
	\end{lemma}
	
	\begin{proof}
		Set
		\[
		F=g_2P_{\sigma^r\omega}^{i-r}g_1,
		\qquad
		B_\nu=1+\|g_\nu\|_{\mathrm{Lip}},
		\quad \nu=1,2.
		\]
		We first decompose $Fh_{\sigma^i\omega}$ into four non-increasing functions.
		For any real Lipschitz function $g$, define
		\[
		\alpha_g(x)
		=
		g(x)+\|g\|_\infty+\operatorname{Lip}(g)(1-x),
		\qquad
		\beta_g(x)
		=
		\|g\|_\infty+\operatorname{Lip}(g)(1-x).
		\]
		Then $g=\alpha_g-\beta_g$, and both functions are non-negative.
		For $x<y$, the Lipschitz inequality gives
		\[
		\alpha_g(y)-\alpha_g(x)
		=
		g(y)-g(x)-\operatorname{Lip}(g)(y-x)
		\le0.
		\]
		The function $\beta_g$ is also non-increasing. Furthermore,
		\[
		0\le\alpha_g,\beta_g\le2\|g\|_{\mathrm{Lip}}.
		\]
		
		Define
		\[
		U=
		\mathcal L_{\sigma^r\omega}^{i-r}
		(\alpha_{g_1}h_{\sigma^r\omega}),
		\qquad
		V=
		\mathcal L_{\sigma^r\omega}^{i-r}
		(\beta_{g_1}h_{\sigma^r\omega}).
		\]
		Products of non-negative non-increasing functions are non-increasing,
		and the transfer operators preserve this property. Hence
		$U,V$ are non-negative and non-increasing. Since
		\[
		0\le
		\alpha_{g_1}h_{\sigma^r\omega},
		\beta_{g_1}h_{\sigma^r\omega}
		\le2B_1h_{\sigma^r\omega},
		\]
		it follows that
		\begin{equation}
			\label{eq:asip-domination}
			0\le U,V
			\le
			2B_1\mathcal L_{\sigma^r\omega}^{i-r}h_{\sigma^r\omega}
			=
			2B_1h_{\sigma^i\omega}.
		\end{equation}
		Thus their heights are controlled at time $i$, independently
		of the length $i-r$.
		
		By \eqref{eq:asip-normalized-iterate},
		\begin{align*}
			Fh_{\sigma^i\omega}&= g_2\mathcal L_{\sigma^r\omega}^{i-r}(g_1h_{\sigma^r\omega}) = g_2\mathcal L_{\sigma^r\omega}^{i-r}\Big(\alpha_{g_1}h_{\sigma^r\omega}- \beta_{g_1}h_{\sigma^r\omega}\Big)\\
			&=
			g_2(U-V)
			=
			\alpha_{g_2}U-\alpha_{g_2}V
			-\beta_{g_2}U+\beta_{g_2}V.
		\end{align*}
		Each of the four terms on the right is non-negative and
		non-increasing, with height at most
		$4B_1B_2H(\sigma^i\omega)$. Therefore,
		\[
		\|Fh_{\sigma^i\omega}\|_{\mathrm{BV}}
		\le32B_1B_2H(\sigma^i\omega).
		\]
		The $L^\infty$ contraction property of $P_\omega$ also yields
		\[
		\|F\|_\infty\le B_1B_2,
		\qquad
		|\mu_{\sigma^i\omega}(F)|\le B_1B_2.
		\]
		Consequently, the function
		$
		[F]_{\sigma^i\omega}h_{\sigma^i\omega}
		$
		has zero Lebesgue integral and satisfies
		\begin{equation}
			\label{eq:asip-centered-bv}
			\|[F]_{\sigma^i\omega}h_{\sigma^i\omega}\|_{\mathrm{BV}}
			\le34B_1B_2H(\sigma^i\omega).
		\end{equation}
		
		Now we use the memory loss property to bound $\|P_{\sigma^i\omega}^{j-i}[F]_{\sigma^i\omega}\|_{L^s(\mu_{\sigma^j\omega})}$ for any  $s\ge1$.
		Put $n=j-i$. Using
		\eqref{eq:asip-normalized-l1},
		\eqref{eq:asip-memory-input}, and
		\eqref{eq:asip-centered-bv}, we obtain
		\begin{align*}
			\|P_{\sigma^i\omega}^n[F]_{\sigma^i\omega}\|_{L^1(\mu_{\sigma^j\omega})}
			&= \|\mathcal L_{\sigma^i\omega}^n ([F]_{\sigma^i\omega} h_{\sigma^i\omega})\|_{L^1(m)} \\
			&\le34B_1B_2H(\sigma^i\omega)Z(\sigma^i\omega)e^{-cn}.
		\end{align*}
		On the other hand,
		\[
		\|P_{\sigma^i\omega}^n[F]_{\sigma^i\omega}\|_{L^\infty(\mu_{\sigma^j\omega}) }
		\le2B_1B_2.
		\]
		Define
		\[
		D_n(\omega)
		=
		\min\{1,17H(\omega)Z(\omega)e^{-cn}\}.
		\]
		Therefore, for $1\le s<\infty$,
		\begin{equation}
			\label{eq:asip-interpolation}
			\|P_{\sigma^i\omega}^n[F]_{\sigma^i\omega}\|_{L^s(\mu_{\sigma^j\omega})}
			\le2B_1B_2D_n^{1/s}(\sigma^i\omega).
		\end{equation}
		
		Recall that the random variable $H$ depends only on past coordinates $(\omega_j)_{j<0}$,
		whereas $Z$ depends only on present and future coordinates $(\omega_j)_{j\geq 0}$.
		Thus they are independent. By
		\eqref{eq:asip-sample-height} and $Z\in L^1(\mathbb P)$, there exists a constant $C_0$ such that
		\begin{equation}
			\label{eq:asip-D-mean}
			0\le D_n\le1,
			\qquad
			\mathbb ED_n
			\le17e^{-cn}\mathbb E(HZ)
			=
			17e^{-cn}\mathbb EH\,\mathbb EZ
			\le C_0e^{-cn}.
		\end{equation}
		
		Put $t=\max\{p,s\}$, and  choose
		\[
		\gamma_{s,p}=\frac{c}{2t},
		\]
		where $c$ is  the exponent in the quenched $BV$ memory-loss
		estimate \eqref{eq:asip-memory-input}.
		Since  $0\le D_n\le1$,
		\[
		\|D_n^{1/s}\|_{L^t(\mathbb P)}
		=
		(\mathbb ED_n^{t/s})^{1/t}
		\le
		(\mathbb ED_n)^{1/t}
		\le C_0^{1/t}e^{-cn/t}.
		\]
		Define
		\[
		K_{s,p}(\omega)
		=
		2\sum_{n=0}^{\infty}
		e^{\gamma_{s,p} n}D_n(\omega)^{1/s}.
		\]
		Minkowski's inequality, first for finite partial sums and then
		by monotone convergence, gives
		\[
		\|K_{s,p}\|_{L^t}
		\le
		2C_0^{1/t}
		\sum_{n=0}^{\infty}e^{(\gamma_{s,p}-c/t)n}
		=
		\frac{2C_0^{1/t}}{1-e^{-c/(2t)}}<\infty.
		\]
		As $\mathbb P$ is a probability measure and $p\le t$,
		this also proves $K_{s,p}\in L^p(\mathbb P)$.
		
		Let $A_{s,p}$ be a measurable set of full
		$\mathbb P$-measure on which Lemma \ref{lem:asip-sample-densities} hold, and the uniform $BV$ memory-loss estimate \eqref{eq:asip-memory-input}
		used above hold, and on which $K_{s,p}$ is finite. Such a set
		exists because $K_{s,p}\in L^p(\mathbb P)$ and the preceding
		estimates hold on common full-measure sets. Define
		\[
		\Omega_{s,p}^*
		=
		\bigcap_{k\in\mathbb Z}\sigma^{-k}A_{s,p}
		=
		\{\omega:\sigma^k\omega\in A_{s,p}
		\text{ for every }k\in\mathbb Z\}.
		\]
		Since $\sigma$ preserves $\mathbb P$, each set
		$\sigma^{-k}A_{s,p}$ has measure one. The intersection is
		countable, so $\mathbb P(\Omega_{s,p}^*)=1$. Moreover,
		$\Omega_{s,p}^*$ is invariant under $\sigma$ and $\sigma^{-1}$.
		Consequently, for every $\omega\in\Omega_{s,p}^*$, all the
		properties used in the proof are available at every shifted
		environment $\sigma^k\omega$, and
		\[
		K_{s,p}(\sigma^k\omega)<\infty
		\qquad\text{for every }k\in\mathbb Z.
		\]
		Thus the argument applies to every choice of integer time
		indices $r\le i\le j$ without changing the exceptional set.
		
		The same exceptional set also works for all Lipschitz
		observables $g_1,g_2$. Indeed, the $BV$ memory-loss estimate
		holds simultaneously for all zero-integral $BV$ functions,
		and the subsequent decompositions and inequalities are
		deterministic once the $\omega$ is fixed. Furthermore,
		$H$, $Z$, $D_n$, and $K_{s,p}$ are constructed independently of
		the observables. This proves that
		\eqref{eq:asip-two-observable} holds on $\Omega_{s,p}^*$
		simultaneously for all the stated observables and time
		indices, with the same $\gamma_{s,p}$ and the same $K_{s,p}$.

	\end{proof}
	
	\begin{theorem}
		\label{thm:quenched-asip}
		Assume the hypotheses of Theorem~\ref{thm:main}, and let
		$\mu_\omega$ be the sample measures constructed in
		Lemma~\ref{lem:asip-sample-densities}. For a  Lipschitz
		function $\psi:[0,1]\to\mathbb R$, define
		\begin{align}\label{eq:sum-var}
			\varphi_\omega=\psi-\mu_\omega(\psi),
			\qquad
			S_n^\omega
			=
			\sum_{j=0}^{n-1}
			\varphi_{\sigma^j\omega}\circ T_\omega^j,
			\qquad
			v_n(\omega)
			=
			\int(S_n^\omega)^2\,d\mu_\omega.
		\end{align}
		There exists a deterministic $\Sigma^2\ge0$ such that
		\[
		\lim_{n\to\infty}\frac{v_n(\omega)}n=\Sigma^2
		\qquad\text{for $\mathbb P$-almost every $\omega$.}
		\]
		If $\Sigma^2>0$, there exists $\delta>0$ such that, for
		$\mathbb P$-almost every $\omega$, a standard Brownian motion
		$B^\omega$ can be defined on an extension of
		$([0,1],\mathcal B([0,1]),\mu_\omega)$ for which
		\[
		S_n^\omega-B^\omega_{v_n(\omega)}
		=
		O(n^{1/2-\delta})
		\qquad\text{almost surely.}
		\]
	\end{theorem}
	
	\begin{proof}
		We verify the remaining hypotheses of
		\cite[Theorem~9]{DHL2026}.
		By construction,
		\[
		\mu_\omega(\varphi_\omega)=0.
		\]
		Moreover, since $\mu_\omega$ is a probability measure,
		\[
		|\mu_\omega(\psi)|\le\|\psi\|_\infty,
		\qquad
		\|\varphi_\omega\|_\infty\le2\|\psi\|_\infty.
		\]
		Subtracting a  constant function does not change the
		Lipschitz seminorm, so
		\[
		\operatorname{Lip}(\varphi_\omega)
		=
		\operatorname{Lip}(\psi).
		\]
		Therefore
		\begin{equation}
			\label{eq:asip-observable-bound}
			\|\varphi_\omega\|_{\mathrm{Lip}}
			\le
			2\|\psi\|_\infty+\operatorname{Lip}(\psi)
			\le2\|\psi\|_{\mathrm{Lip}}
		\end{equation}
		uniformly in $\omega$. Thus the random variable
		$\omega\mapsto\|\varphi_\omega\|_{\mathrm{Lip}}$ belongs to
		$L^\infty(\mathbb P)$. In the notation of \cite[Theorem~9]{DHL2026},
		the observable integrability exponent may consequently be
		taken to be $r=\infty$.

		Apply Lemma~\ref{lem:asip-two-observable} with $p=12$, $s \in \{2,3\}$,
		and write
		\[
		K= \max\{  K_{2, 12}, K_{3, 12}   \} ,
		\qquad
		\gamma= \min\{ \gamma_{2,12}, \gamma_{3,12} \}.
		\]
		The lemma verifies Assumption~1 of \cite{DHL2026} with
		\[
		p=12,\qquad
		I=\{2,3\},
		\qquad
		a_{2,n}=a_{3,n}=e^{-\gamma n}.
		\]
		Choose $\widetilde s=3/2$, so that
		$I=\{2,2\widetilde s\}$. The first two numerical conditions in
		\cite[Theorem~9]{DHL2026} become
		\[
		\frac1p+\frac1{r}
		=
		\frac1{12}<\frac14
		\]
		and
		\[
		2>\widetilde s=\frac32
		>
		\frac1{1-2/p-2/r}
		=
		\frac65.
		\]
		Since exponential decay implies
		$a_{2,n}=O(n^{-2})$, we may take $b=2$. The next condition is
		then satisfied because
		\[
		\frac{2/p+4/r}{b-1}
		=
		\frac16
		<
		\frac56
		=
		1-\frac4{r}-\frac2p.
		\]
		Finally,
		\[
		\sum_{n=1}^{\infty}a_{2\widetilde s,n}
		=
		\sum_{n=1}^{\infty}e^{-\gamma n}<\infty.
		\]
		The variance limit and the asserted ASIP therefore follow
		from \cite[Theorem~9]{DHL2026}.
	\end{proof}
	
	\begin{theorem}
		Under the assumptions of Theorem
		\ref{thm:quenched-asip}, for any $\varepsilon > 0$ and for $\mathbb P$-almost every $\omega$,
		\begin{align*}
			\sup_{ t \in \mathbb{R} } |  \mu_\omega( S_n^\omega \le  t \sqrt{v_n(\omega)}  )  - \Phi(t)  | = O( n^{-1/5 + \varepsilon} ).
		\end{align*}
		Here, $S_n^\omega$ and $v_n(\omega)$ are defined as in \eqref{eq:sum-var}, and $\Phi(t)$ denotes the standard normal
		distribution function.
	\end{theorem}
	
	\begin{proof} By Lemma~\ref{lem:asip-two-observable}, Assumption 1 of \cite{DHL2026} holds with
		$$
		\mathcal{I} = \{ 2, 2 \tilde{s} \}, \quad \frac{1}{1 - 2/p} < \tilde{s} < 2, \quad a_{2,n} = a_{ 2 \tilde{s} , n  }
		= e^{  - \gamma_{ \tilde{s}, p  } n  },
		$$
		where $p \ge 1$ can be taken arbitrarily large. Proposition 19 in \cite{DHL2026} implies that, for any $\varepsilon > 0$,
		\begin{align*}
			\biggl\Vert \sum_{j=0}^{n-1} ( H_{\omega,j} \circ T_\omega^j )^2 - \sum_{j=0}^{n-1} \mathbb{E}_{ \sigma^j \omega }
			[ (H_{\omega,j})^2 ] \biggl\Vert_{L^{ \tilde{s} }(\mu_\omega)}^{ \tilde{s} }
			= O( n^{  \tilde{s} /2 + \varepsilon } + n^{ - \tilde{s} + 2  + \varepsilon }  ),
		\end{align*}
		where $H_{\omega,j}$ are random variables on $( [0,1], \mu_\omega )$ defined in \cite[equation (14)]{DHL2026}.
		Consequently, taking $\tilde{s}$ sufficiently close to $2$, we see that
		the condition in \cite[Proposition 18]{DHL2026} holds with $s = \tilde{s}$ and $A = 1$.
		The desired estimate then follows from \cite[Proposition 18]{DHL2026}.
	\end{proof}
	
	\section*{Acknowledgements}
	H. Cui was supported by the National Natural Science Foundation of China (Nos. 12271418) and the Key project of Hubei Provincial Department of Education (Grant No. D20231103).
	J. Lepp\"anen was supported by JSPS via the project LEADER. He would like to thank the Wuhan University of
	Science and Technology for hospitality over the course of this work.
	\vspace{0.3cm}

	
	
	


\begin{thebibliography}{99}
		
		\bibitem{BBR19} W. Bahsoun, C. Bose and M. Ruziboev,
		Quenched decay of correlations for slowly mixing systems,
		\emph{Trans. Amer. Math. Soc.} \textbf{372} (2019), no. 9, 6547--658
		
		\bibitem{B99} J. Buzzi, Exponential decay of correlations for random Lasota--Yorke maps,
		\emph{Commun. Math. Phys.} \textbf{208} (1999), no. 1, 25--54
		
		\bibitem{Cui2021}
		H.~Cui,
		\newblock Decay of correlations and memory loss for Lasota--Yorke convex maps,
		\newblock \emph{Dyn. Syst.} \textbf{36} (2021), no.~3, 404--415
		
		\bibitem{CS24} H. Cui and J. Sui, Large deviation and central limit theorems for sequential and random Lasota--Yorke convex maps,
		\emph{Fundam. Math.} \textbf{267} (2024), 267--289
		
		\bibitem{CW25} H. Cui and C. Wang,
		Maximal large deviations for sequential dynamical systems,
		\emph{J. Stat. Phys.} \textbf{192} (2025), art. 163
		
		
		
		\bibitem{CR07}
		J.-P. Conze and A. Raugi, Limit theorems for sequential expanding dynamical systems on $[0, 1]$,
		\emph{Contemp. Math.} {\bf 430} (2007), 89--121
		
		\bibitem{DS23} D. Dragi\v{c}evi\'{c} and J. Sedro,
		Quenched limit theorems for expanding on average cocycles,
		\emph{Stoch. Dyn.} \textbf{23} (2023), no. 6, 2350046
		
		\bibitem{DH21} D. Dragi\v{c}evi\'{c} and Y. Hafouta,
		Almost sure invariance principle for random distance expanding maps
		with a nonuniform decay of correlations,
		\emph{Lect. Notes Math.} \textbf{2290} (2021), 177--204
		
		\bibitem{DHS23} D. Dragi\v{c}evi\'{c}, Y. Hafouta and J. Sedro,
		A vector-valued almost sure invariance principle for random expanding on average cocycles,
		\emph{J. Stat. Phys.} \textbf{190} (2023), art.~54.
		
		\bibitem{DHL2026}
		D.~Dragi\v{c}evi\'{c}, Y.~Hafouta and J.~Lepp\"{a}nen,
		\newblock Decay of correlations and limit theorems for random intermittent maps,
		\newblock \emph{J. Lond. Math. Soc.} \textbf{114} (2026), no. 3, e70707
		
		\bibitem{GOT13} C. Gupta, W. Ott and A. T\"or\"ok,
		Memory loss for time-dependent piecewise expanding systems in higher dimension,
		\emph{Math. Res. Lett.} {\bf 20} (2013), no. 1, 141--161
		
		\bibitem{H23}
		Y. Hafouta, Spectral methods for limit theorems for random expanding transformations,
		\emph{preprint} (2023), arXiv:2311.12950
		
		\bibitem{Kallenberg2002}
		O. ~Kallenberg, Foundations of Modern Probability, 2nd edition, Springer, 2002
		
		\bibitem{KL21}
		A. Korepanov and J. Lepp\"anen,
		Loss of memory and moment bounds for nonstationary intermittent dynamical systems,
		\emph{Commun. Math. Phys.} \textbf{385} (2021), no. 2, 905--935
		
		\bibitem{OSY09} W. Ott, M. Stenlund and L.-S. Young,
		Memory loss for time-dependent dynamical systems,
		\emph{Math. Res. Lett.} \textbf{16} (2009), no. 3,
		463--475
		
		\bibitem{StenlundSulku2014}
		M.~Stenlund and H.~Sulku,
		\newblock A coupling approach to random circle maps expanding on the average,
		\newblock \emph{Stoch. Dyn.} \textbf{14} (2014), no.~4, 1450008
		
		\bibitem{SYZ13} M. Stenlund, L.-S. Young, H. Zhang,
		Dispersing billiards with moving scatterers.
		\emph{Commun. Math. Phys.} \textbf{322} (2013), no. 3, 909--955
		
	\end{thebibliography}
\end{document}